\documentclass[reqno]{amsart}

\usepackage{amssymb,amsthm,fancyheadings,amsmath,dsfont}

\newtheorem{defi}{Definition}[section]
\newtheorem{thm}[defi]{Theorem}

\newtheorem{cor}[defi]{Corollary}

\newtheoremstyle{normal}
    {}              
    {}              
    {\normalfont}       
    {}                  
    {\bfseries}         
    {.}                 
    {.5em}              
    {}                  
\theoremstyle{normal}

\newtheorem{rem}[defi]{Remark}

\newcommand{\bpr}{\begin{proof}[Proof]}  
\newcommand{\epr}{\end{proof}}
\newcommand{\beq}{\begin{equation}}
\newcommand{\eeq}{\end{equation}}
\newcommand{\bce}{\begin{center}}
\newcommand{\ece}{\end{center}}
\newcommand{\be}{\begin{enumerate}}  
\newcommand{\ee}{\end{enumerate}}

\newcommand{\compemb}{\lhook\hspace{-0.151cm}{-}\hspace{-0.3cm}\hookrightarrow}

\DeclareMathOperator*{\dist}{dist}
\DeclareMathOperator*{\diver}{div}

\def\pa{\partial}
\def\om{\omega}
\def\si{\sigma}

\def\ep{\varepsilon}
\def\de{\delta}
\def\ph{\varphi}

\def\ga{\gamma}

\def\De{\Delta}

\def\Om{\Omega}

\def\R{\mathbb R}
\def\C{\mathbb C}

\def\N{\mathbb N}

\def\B{\mathbb B}
\def\E{\mathbb E}

\def\calA{\mathcal A}
\def\calB{\mathcal B}

\def\calT{\mathcal T}
\def\calE{\mathcal E}
\def\calR{\mathcal R}

\def\calM{\mathcal M}

\def\calP{\mathcal P}

\def\calU{\mathcal U}
\def\calK{\mathcal K}
\def\calN{\mathcal N}
\def\calV{\mathcal V}

\numberwithin{equation}{section}

\title[Quasilinear Parabolic Evolution Equations in weighted $L_p$ spaces]
{On quasilinear Parabolic Evolution Equations \\ in weighted $L_p$-spaces}

\author[Matthias K\"{o}hne]{Matthias K\"{o}hne}
\address{Center for Computational Engineering Science, RWTH Aachen University, Pauwelsstr.~19, 52074 Aachen, Germany}
\email{koehne@mathcces.rwth-aachen.de}

\author[Jan Pr\"uss]{Jan Pr\"uss}
\address{Institut f\"ur Mathematik, Martin-Luther-Universit\"at Halle-Wittenberg,
Theodor-Lieser-Str. 5, 06120 Halle, Germany}
\email{jan.pruess@mathematik.uni-halle.de}

\author[Mathias Wilke]{Mathias Wilke}
\address{Institut f\"ur Mathematik, Martin-Luther-Universit\"at Halle-Wittenberg,
Theodor-Lieser-Str. 5, 06120 Halle, Germany}
\email{mathias.wilke@mathematik.uni-halle.de (corresponding author)}

\subjclass[2000]{}

\date{\today}


\begin{document}

\maketitle

\begin{abstract}
In this paper we develop a geometric theory for quasilinear parabolic problems in weighted $L_p$-spaces. We prove existence and uniqueness of solutions as well as the continuous dependence on the initial data. Moreover, we make use of a regularization effect for quasilinear parabolic equations to study the $\om$-limit sets and the long-time behaviour of the solutions. These techniques are applied to a free boundary value problem. The results in this paper are mainly based on maximal regularity tools in (weighted) $L_p$-spaces.
\end{abstract}

\section{Introduction}

\noindent
In this paper we consider abstract quasilinear parabolic problems of the form
    \begin{equation}\label{Int1}
    \dot{u}+A(u)u=F(u),\ t>0,\quad u(0)=u_0,
    \end{equation}
where $(A,F):V_\mu\to \calB(X_1,X_0)\times X_0$ and $u_0\in V_\mu$. The spaces $X_1,X_0$ are Banach spaces such that $X_1\hookrightarrow X_0$ with dense embedding and $V_\mu$ is an open subset of the real interpolation space
    $$X_{\gamma,\mu}:=(X_0,X_1)_{\mu-1/p,p},\quad \mu\in(1/p,1].$$
By $\calB(X_1,X_0)$ we denote the set of all bounded linear operators from $X_1$ to $X_0$. For $p\in (1,\infty)$, let $L_{p,\mu}(J;X)$ denote the vector-valued weighted $L_p$-space
    \begin{equation}\label{Int2}
    L_{p,\mu}(J;X):=\{u:J\to X_0:t^{1-\mu}u\in L_p(J;X)\},
    \end{equation}
where $X$ is a Banach space, $\mu\in (1/p,1]$ and $J=[0,T]$, $T>0$. In this paper we are interested in solutions $u(t)$ of \eqref{Int1} having \emph{maximal} $L_{p,\mu}$-\emph{regularity}, i.e.
    $$u\in H_{p,\mu}^1(J;X_0)\cap L_{p,\mu}(J;X_1),$$
with $H_{p,\mu}^1(J;X_0)$ being defined as
    $$H_{p,\mu}^1(J;X_0):=\{u\in L_{p,\mu}(J;X_0)\cap W_{1}^1(J;X_0):\dot{u}\in L_{p,\mu}(J;X_0)\},$$
and $H_{p,\mu}^1(J;X_0)$ is supplied with the norm
    $$||u||_{H_{p,\mu}^{1}}:=||u||_{L_{p,\mu}}+||\dot{u}||_{L_{p,\mu}},$$
which turns it into a Banach space. In our approach it is crucial to know that the operator $A_0:=A(u_0)$ has the property of maximal $L_{p,\mu}$-regularity, for short $A_0\in\calM\calR_{p,\mu}(X_1,X_0)$. To be precise, this means that for each $f\in L_{p,\mu}(\R_+;X_0)$ there exists a unique solution
    $$u\in H_{p,\mu}^1(\R_+;X_0)\cap L_{p,\mu}(\R_+;X_1)$$
of the problem
    $$\dot{u}+A_0u=f,\ t>0,\quad u(0)=0.$$
Thanks to \cite[Theorem 2.4]{PrSi04} the characterization
    $$A_0\in \calM\calR_{p,\mu}(X_1,X_0)\Leftrightarrow A_0\in \calM\calR_p(X_1,X_0)$$
for a closed linear operator $A_0$ in $X_0$ holds true, provided $\mu\in (1/p,1]$, $p\in (1,\infty)$. Here we use the notation $A_0\in \calM\calR_p(X_1,X_0)$ for the 'classical' case $\mu=1$. This characterization is very useful, since there are many results available which ensure $A_0\in \calM\calR_p(X_1,X_0)$, see e.g.\ \cite{DHP1}. Concerning nontrivial initial data, it was shown in \cite[Theorem 3.2]{PrSi04} that if $A_0\in \calM\calR_p(X_1,X_0)$, then the initial value problem
    $$\dot{u}+A_0u=f,\ t>0,\quad u(0)=u_0.$$
has a unique solution $u$ with maximal $L_{p,\mu}$-regularity if and only if $f\in L_{p,\mu}(\R_+;X_0)$ and $u_0\in X_{\gamma,\mu}$, which is the natural phase space in this functional analytic setting.

The choice of the weighted $L_p$-spaces has a big advantage. To see this, observe that for each fixed $\de\in (0,T)$ the embedding
    $$H_{p,\mu}^1(0,T;X_0)\cap L_{p,\mu}(0,T;X_1)\hookrightarrow H_{p}^1(\de,T;X_0)\cap L_{p}(\de,T;X_1)$$
is true. This shows that if we start with an initial value in the larger space $X_{\gamma,\mu},\ \mu\in (1/p,1)$, compared to the classical case $\mu=1$, the solution regularizes instantaneously, since $\de>0$ may be arbitrarily small. Note that this regularization effect can not be obtained in the usual setting of maximal $L_p$-regularity, i.e.\ if $\mu=1$. We use this property to study the long-time behaviour of the solutions of \eqref{Int1}, in particular their $\om$-limit sets. To our knowledge, so far, there do not exist results on well-posedness of \eqref{Int1} and its consequences in weighted $L_p$-spaces of the form \eqref{Int2}.

This paper is organized as follows. In Section 2 we show that the initial value problem \eqref{Int1} has maximal $L_{p,\mu}$-regularity, if $A(u_0)\in\calM\calR_{p}(X_1,X_0)$ and if $(A,F)$ are Lipschitz continuous. Furthermore we show that the solutions to \eqref{Int1} depend continuously on the initial data. These results extend those of Cl\'{e}ment \& Li \cite{CleLi93} and Pr\"{u}ss \cite{JanBari} who considered unweighted $L_p$-spaces, i.e.\ the case $\mu=1$. \\
In Section 3 we prove that bounded orbits in $X_{\gamma}:=X_{\gamma,1}$ are already relatively compact in $X_{\gamma}$, provided $X_\gamma$ is compactly embedded in $X_{\gamma,\mu}$, $\mu\in (1/p,1).$
In particular this yields global existence of solutions which are bounded in $X_\gamma$. By means of the variation of parameters formula, this is easy to prove for semilinear equations, where $A(u)\equiv A_0$, but in the quasilinear case it is by no means obvious. For this purpose we make use of the regularization effect as well as of the continuous dependence of the solutions on the initial data. At the end of Section 3 we apply this result to a class of second order quasilinear parabolic initial boundary value problems.\\
Section 4 is devoted to the long-time behaviour of solutions of \eqref{Int1}. By relative compactness of the orbits, the $\om$-limit set $\om(u_0)\subset X_{\gamma}$ of the solution $u(t)$ to \eqref{Int1} is nonempty, compact, connected and a global attractor for the solution $u(t)$. Assuming the existence of a strict Ljapunov functional, we have furthermore $\om(u_0)\subset\calE$, where $\calE$ denotes the set of equilibria of $\eqref{Int1}$, i.e.\ the set of all solutions of \eqref{Int1} which are constant in $t$. If we postulate that there exists $u_*\in\om(u_0)$ which is \emph{normally hyperbolic} (see Theorem \ref{thmPSZ} for the notion of normal hyperbolicity) , it follows that $u(t)$ converges at an exponential rate to $u_*$ in $X_\gamma$, provided $(A,F)$ are continuously differentiable. In this way we extend the local convergence result \cite[Theorem 6.1]{PSZ} to a global one, i.e.\ there is no need to choose the initial value sufficiently close to $u_*$ in $X_\gamma$.\\
Finally, in Section 5, we show that the techniques of Section 3 \& 4 can also be applied to problems with moving boundary. To be precise, we study global existence and long-time behaviour of solutions to the two-phase Mullins-Sekerka problem. For the sake of readability and completeness we also provide some facts from differential geometry, which are essential for our considerations.

There is a vast literature concerning existence and uniqueness of solutions to quasilinear parabolic problems of the form \eqref{Int1} in different functional analytic settings, see \cite{Ama88,Ama89,Ama90,Ama93,Ama05,Ang90,CleLi93,CleSi01,Esch94,LSU,Lun95,Sim94,Sim95,Yag91}; this is just a selection.  Basically there are two approaches to establish well-posedness of \eqref{Int1}. One makes use of the theory of parabolic evolution operators, see e.g.\ \cite{Ama88,Ama93}. Another approach uses maximal regularity tools which have for instance been applied in \cite{Ama05,Ang90,CleLi93,CleSi01,JanBari,Sim94}. The method of maximal regularity has the advantage that it provides a natural analytic setting for the semiflow, which is induced by \eqref{Int1}. A theory based on function spaces with weights has been used in \cite{Ang90} in order to treat functions with a singularity at $t=0$. This approach has been further developed in the papers \cite{Sim94} and \cite{CleSi01}, which are based on maximal regularity in continuous interpolation spaces. In particular, the authors in \cite{CleSi01} consider
    \begin{multline*}
    BUC_{\mu}([0,T];X):=\{u\in C((0,T];X):t^{1-\mu}u\in BUC((0,T];X),\\
    \lim_{t\to 0+}t^{1-\mu}|u(t)|_X=0\},\quad \mu\in (0,1].
    \end{multline*}
as a basic space, instead of \eqref{Int2}. Concerning the long-time behaviour of solutions, we refer e.g.\ to \cite{BrHuLu00,LaPrSchn06,LaPrSchn08,Lun95,PSZ,Sim95}. It is one aim of this paper to extend the local convergence result \cite[Theorem 6.1]{PSZ} to a global one. At this point we want to mention the paper \cite{ChFaSch09} where the authors prove a Lojasiewicz inequality for the Willmore flow, a problem for moving hypersurfaces. They apply this inequality to exclude compact blowups for the Willmore flow.

\textbf{Notations.} Let $p\in (1,\infty)$, $T\in (0,\infty)$ and $\mu\in (1/p,1]$. If $X_0$ and $X_1$ are Banach spaces with dense embedding $X_1\hookrightarrow X_0$, we define
    $$\E_{1,\mu}(0,T):=H_{p,\mu}^1(0,T;X_0)\cap L_{p,\mu}(0,T;X_1),$$
    $$\E_{0,\mu}(0,T):=L_{p,\mu}(0,T;X_0),$$
and
    $$X_{\gamma,\mu}:=(X_0,X_1)_{\mu-1/p,p},$$
where $(X_0,X_1)_{\mu-1/p,p}$ is the real interpolation space of order $\mu-1/p$ and exponent $p$. Furthermore we denote by $||\cdot||_{\infty,X_{\gamma,\mu}}$ the norm in $BC([0,T];X_{\gamma,\mu})$. In the 'classical' case $\mu=1$ we simply use the notation $\E_0$, $\E_1$ and $X_\gamma$ instead of $\E_{1,1}$, $\E_{0,1}$ and $X_{\gamma,1}$. We write $X_1{\compemb} X_0$ if $X_1$ is compactly embedded in $X_0$. If $M_1$ and $M_2$ are metric spaces and $F:M_1\to M_2$, then $F\in C^{1-}(M_1;M_2)$ means that $F$ is locally Lipschitz.

\section{Local Well-Posedness}\label{LWPsec}

The aim of this section is to solve the quasilinear evolution equation
    \begin{equation}\label{LWP1}
    \dot{u}+A(u)u=F(u),\ t>0,\quad u(0)=u_1,
    \end{equation}
under the assumption that there exist two Banach spaces $X_0,X_1$, with dense embedding $X_1\hookrightarrow X_0$ such that the nonlinear mappings $(A,F)$ satisfy
    \begin{equation}\label{LWP2}
    (A,F)\in C^{1-}(V_\mu;\calB(X_1,X_0)\times X_0),
    \end{equation}
where $V_\mu\subset (X_0,X_1)_{\mu-1/p,p}=:X_{\gamma,\mu}$ is open and nonempty for some $\mu\in (1/p,1]$. The main result of this section reads as follows.
\begin{thm}\label{LWPthm}
Let $p\in (1,\infty)$, $u_0\in V_\mu$ be given and suppose that $(A,F)$ satisfy \eqref{LWP2} for some $\mu\in (1/p,1]$. Assume in addition that $A(u_0)\in\calM\calR_{p}(X_1,X_0)$. Then there exist $T=T(u_0)>0$ and $\ep=\ep(u_0)>0$, such that $\bar{B}_{\ep}^{X_{\gamma,\mu}}(u_0)\subset V_\mu$ and such that problem \eqref{LWP1} has a unique solution
    $$u(\cdot,u_1)\in H_{p,\mu}^1(0,T;X_0)\cap L_{p,\mu}(0,T;X_1)\cap C([0,T];V_\mu),$$
on $[0,T]$, for any initial value $u_1\in \bar{B}_\ep^{X_{\gamma,\mu}}(u_0)$. Furthermore there exists a constant $c=c(u_0)>0$ such that for all $u_1,u_2\in \bar{B}_\ep^{X_{\gamma,\mu}}(u_0)$ the estimate
    $$||u(\cdot,u_1)-u(\cdot,u_2)||_{\E_{1,\mu}(0,T)}\le c|u_1-u_2|_{X_{\gamma,\mu}}$$
is valid.
\end{thm}
\bpr
Since $u_0\in V_\mu$ and by \eqref{LWP2}, there exists $\ep_0>0$ and a constant $L>0$ such that $\bar{B}_{\ep_0}^{X_{\gamma,\mu}}(u_0)\subset V_\mu$ and
    \begin{equation}\label{LWP3}
    |A(w_1)v-A(w_2)v|_{X_0}\le L|w_1-w_2|_{X_{\gamma,\mu}}|v|_{X_1},
    \end{equation}
as well as
    \begin{equation}\label{LWP4}
    |F(w_1)-F(w_2)|_{X_0}\le L|w_1-w_2|_{X_{\gamma,\mu}},
    \end{equation}
hold for all $w_1,w_2\in \bar{B}_{\ep_0}^{X_{\gamma,\mu}}(u_0)$, $v\in X_1$. By the results of the previous section we may introduce a reference function $u_0^*\in \E_{1,\mu}(0,T)$ as the solution of the linear problem
    $$\dot{w}+A(u_0)w=0,\quad w(0)=u_0.$$
Define a ball $\B_r\subset \E_{1,\mu}(0,T)$ by
    $$\B_{r,T,u_1}:=\{v\in \E_{1,\mu}(0,T):v|_{t=0}=u_1\ \mbox{and}\ ||v-u_0^*||_1\le r\},\quad 0<r\le1.$$
Let $u_1\in \bar{B}_{\ep}^{X_{\gamma,\mu}}(u_0)$ with $\ep\in (0,\ep_0]$. We will show that for all $v\in\B_{r,T,u_1}$ it holds that $v(t)\in \bar{B}_{\ep_0}^{X_{\gamma,\mu}}(u_0)$ for all $t\in [0,T]$, provided that $r,T,\ep>0$ are sufficiently small. For this purpose we define $u_1^*\in\E_{1,\mu}(0,T)$ as the unique solution of
    $$\dot{w}+A(u_0)w=0,\quad w(0)=u_1.$$
Given $v\in \B_{r,T,u_1}$ we estimate as follows.
    \begin{align}\label{LWP4a}
      ||v-u_0||_{\infty,X_{\gamma,\mu}}\le||v-u_1^*||_{\infty,X_{\gamma,\mu}}+||u_1^*-u_0^*||_{\infty,X_{\gamma,\mu}}
        +||u_0^*-u_0||_{\infty,X_{\gamma,\mu}}.
    \end{align}
Since $u_0$ is fixed, there exists $T_0=T_0(u_0)>0$ such that $\sup_{t\in[0,T_0]}|u_0^*(t)-u_0|_{X_{\gamma,\mu}}\le\ep_0/3$. Observe that $(v-u_1^*)|_{t=0}=0$, hence
    $$||v-u_1^*||_{\infty,X_{\gamma,\mu}}\le C_1||v-u_1^*||_{\E_{1,\mu}(0,T)}$$
and the constant $C_1>0$ does not depend on $T$. Therefore
    \begin{align*}
    ||v-u_1^*||_{\infty,X_{\gamma,\mu}}&\le C_1||v-u_1^*||_{\E_{1,\mu}(0,T)}\le C_1(||v-u_0^*||_{\E_{1,\mu}(0,T)}+||u_0^*-u_1^*||_{\E_{1,\mu}(0,T)})\\
        &\le C_1(r+||u_0^*-u_1^*||_{\E_{1,\mu}(0,T)}),
    \end{align*}
and \eqref{LWP4a} yields the estimate
    \begin{multline*}
    ||v-u_0||_{\infty,X_{\gamma,\mu}}\\
        \le C_1(r+||u_0^*-u_1^*||_{\E_{1,\mu}(0,T)})+||u_0^*-u_1^*||_{\infty,X_{\gamma,\mu}}+||u_0^*-u_0||_{\infty,X_{\gamma,\mu}}.
    \end{multline*}
Since by assumption the semigroup $e^{-A(u_0)t}$ is exponentially stable it follows that
    \begin{equation}\label{LWP4b}
    ||u_0^*-u_1^*||_{\infty,X_{\gamma,\mu}}+C_1||u_0^*-u_1^*||_{\E_{1,\mu}(0,T)}\le C_\gamma|u_0-u_1|_{X_{\gamma,\mu}},
    \end{equation}
with a constant $C_\gamma>0$ which does not depend on $T$. Choosing $\ep\le\ep_0/(3C_\gamma)$ and $r\le\ep_0/(3C_1)$, we finally obtain
    \begin{equation}\label{LWP4c}
    ||v-u_0||_{\infty,X_{\gamma,\mu}}\le C_1r+C_\gamma\ep+||u_0^*-u_0||_{\infty,X_{\gamma,\mu}}\le\ep_0.
    \end{equation}
Throughout the remainder of this proof we will assume that $u_1\in B_{\ep}^{X_{\gamma,\mu}}(u_0)$, $\ep\le\ep_0/(3C_\gamma)$, $T\in[0,T_0]$ and $r\le\ep_0/(3C_1)$. Under these assumptions, we may define a mapping $\calT_{u_1}:\B_{r,T,u_1}\to \E_{1,\mu}(0,T)$ by means of $\calT_{u_1} v=u$, where $u$ is the unique solution of the linear problem
    \begin{equation*}
    \dot{u}+A(u_0)u=F(v)+(A(u_0)-A(v))v,\ t>0,\quad u(0)=u_1.
    \end{equation*}
In order to apply the contraction mapping principle, we have to show $\calT_{u_1}\B_{r,T,u_1}\subset\B_{r,T,u_1}$ and that $\calT_{u_1}$ defines a strict contraction on $\B_{r,T,u_1}$, i.e. there exists $\kappa\in(0,1)$ such that
    $$||\calT_{u_1} v-\calT_{u_1}\bar{v}||_{\E_{1,\mu}(0,T)}\le \kappa||v-\bar{v}||_{\E_{1,\mu}(0,T)},$$
is valid for all $v,\bar{v}\in\B_{r,T,u_1}$. We will first take care about the self-mapping property. Note that for $v\in\B_{r,T,u_1}$ we have
    $$(\calT_{u_1} v)(t)-u_0^*(t)=u_1^*(t)-u_0^*(t)+\left(e^{-A(u_0)\cdot}\ast(F(v)+(A(u_0)-A(v))v)\right)(t).$$
To treat the convolution term, we observe $\left(e^{-A(u_0)\cdot}\ast(F(v)+(A(u_0)-A(v))v)\right)(0)=0$, hence $A(u_0)\in\calM\calR_p(X_1,X_0)$ implies
    \begin{multline*}
    ||e^{-A(u_0)\cdot}\ast(F(v)+(A(u_0)-A(v))v)||_{\E_{1,\mu}(0,T)}\\
        \le C_0||F(v)+(A(u_0)-A(v))v)||_{\E_{0,\mu}(0,T)},
    \end{multline*}
and $C_0>0$ does not depend on $T$. Let us first estimate $(A(u_0)-A(v))v$ in $\E_0(0,T)$. By \eqref{LWP3} and \eqref{LWP4c} we obtain
    \begin{align*}
    ||(A(u_0)-A(v))v||_{\E_{0,\mu}(0,T)}&\le L||v-u_0||_{\infty,X_{\gamma,\mu}}||v||_{\E_{1,\mu}(0,T)}\\
    &\le L||v-u_0||_{\infty,X_{\gamma,\mu}}(r+||u_0^*||_{\E_{1,\mu}(0,T)})\\
    &\le L(C_1r+C_\gamma\ep+||u_0^*-u_0||_{\infty,X_{\gamma,\mu}})(r+||u_0^*||_{\E_{1,\mu}(0,T)}).
    \end{align*}
Furthermore, by \eqref{LWP4} and \eqref{LWP4c}
    \begin{align*}
    ||F(v)||_{\E_{0,\mu}(0,T)}&\le ||F(v)-F(u_0)||_{\E_{0,\mu}(0,T)}+||F(u_0)||_{\E_{0,\mu}(0,T)}\\
    &\le \sigma(T) L||v-u_0||_{\infty,X_{\gamma,\mu}}+||F(u_0)||_{\E_{0,\mu}(0,T)}\\
    &\le \sigma(T) L(C_1r+C_\gamma\ep+||u_0^*-u_0||_{\infty,X_{\gamma,\mu}})+||F(u_0)||_{\E_{0,\mu}(0,T)}\\
    &=\sigma(T)\left[L(C_1r+C_\gamma\ep+||u_0^*-u_0||_{\infty,X_{\gamma,\mu}})+|F(u_0)|_{X_0}\right],
    \end{align*}
with $\sigma(T):=\frac{1}{(1+(1-\mu)p)^{1/p}}T^{1/p+1-\mu}$. Since
    $$||u_0^*-u_0||_{\infty,X_{\gamma,\mu}},||u_0^*||_{\E_{1,\mu}(0,T)}\to 0\ \mbox{as}\ T\to 0_+,$$
this yields
    $$||\calT_{u_1}v-u_0^*||_{\E_{1,\mu}(0,T)}\le ||u_1^*-u_0^*||_{\E_{1,\mu}(0,T)}+r/2,$$
provided $r>0,T>0,\ep>0$ are chosen properly. By \eqref{LWP4b} we obtain in addition
    $$||\calT_{u_1}v-u_0^*||_{\E_{1,\mu}(0,T)}\le (C_\gamma/C_1)|u_1-u_0|_{X_{\ga,\mu}}+r/2\le r/2+r/2=r,$$
with a probably smaller $\ep>0$. This proves the self-mapping property of $\calT_{u_1}$.

Let $u_1,u_2\in \bar{B}_{\ep}^{X_{\ga,\mu}}(u_0)$ be given and let $v_1\in \B_{r,T,u_1}$, $v_2\in \B_{r,T,u_2}$. Then, since $A(u_0)\in\calM\calR_p(X_1,X_0)$, we have
    \begin{multline}\label{LWP5}
    ||\calT_{u_1}v_1-\calT_{u_2}v_2||_{\E_{1,\mu}(0,T)}\le ||e^{-A(u_0)\cdot}(u_1-u_2)||_{\E_{1,\mu}(0,T)}+C_0||F(v_1)-F(v_2)||_{\E_{0,\mu}(0,T)}\\
        +C_0||(A(v_1)-A(u_0))(v_1-v_2)||_{\E_{0,\mu}(0,T)}+C_0||(A(v_1)-A(v_2))v_2||_{\E_{0,\mu}(0,T)}.
    \end{multline}
For the first term on the right hand side we can make use of \eqref{LWP4b} where $u_0$ and $u_0^*$ have to be replaced by $u_2$ and $e^{-A(u_0)t}u_2$, respectively. The second term can be treated as follows. By \eqref{LWP4}, we obtain
    $$||F(v_1)-F(v_2)||_{\E_{0,\mu}(0,T)}\le \sigma(T)L||v_1-v_2||_{\infty,X_{\ga,\mu}}.$$
Moreover, by \eqref{LWP4b} and the trace theorem we have
    \begin{align}\label{LWP6}
    \begin{split}
    ||v_1-v_2||_{\infty,X_{\ga,\mu}}&\le ||v_1-v_2-(e^{-A(u_0)\cdot}(u_1-u_2))||_{\infty,X_{\ga,\mu}}\\
    &\hspace{4cm}+||e^{-A(u_0)\cdot}(u_1-u_2)||_{\infty,X_{\ga,\mu}}\\
    &\le C_1||v_1-v_2-(e^{-A(u_0)\cdot}(u_1-u_2))||_{\E_{1,\mu}(0,T)}+C_\gamma|u_1-u_2|_{X_{\ga,\mu}}\\
    &\le C_1||v_1-v_2||_{\E_{1,\mu}(0,T)}+C_\gamma(1+C_1)|u_1-u_2|_{X_{\ga,\mu}}.
    \end{split}
    \end{align}
This yields
    $$||F(v_1)-F(v_2)||_{\E_{0,\mu}(0,T)}\le \sigma(T)L\left(C_1||v_1-v_2||_{\E_{1,\mu}(0,T)}+C_\gamma(1+C_1)|u_1-u_2|_{X_{\ga,\mu}}\right).$$
For the remaining terms in \eqref{LWP5} we make use of \eqref{LWP3} which results in
    \begin{multline*}
    ||(A(v_1)-A(u_0))(v_1-v_2)||_{\E_{0,\mu}(0,T)}+||(A(v_1)-A(v_2))v_2||_{\E_{0,\mu}(0,T)}\\
        \le L(||v_1-u_0||_{\infty,X_{\ga,\mu}}||v_1-v_2||_{\E_{1,\mu}(0,T)}+||v_1-v_2||_{\infty,X_{\ga,\mu}}||v_2||_{\E_{1,\mu}(0,T)}.
    \end{multline*}
By \eqref{LWP4c}, the term $||v_1-u_0||_{\infty,X_{\ga,\mu}}$ can be made as small as we wish by decreasing $r>0,T>0$ and $\ep>0$. Furthermore we have
    $$||v_2||_{\E_{1,\mu}(0,T)}\le ||v_2-u_0^*||_{\E_{1,\mu}(0,T)}+||u_0^*||_{\E_{1,\mu}(0,T)}\le r+||u_0^*||_{\E_{1,\mu}(0,T)},$$
hence $||v_2||_{\E_{1,\mu}(0,T)}$ is small, provided $r>0$ and $T>0$ are small enough. Lastly, the term $||v_1-v_2||_{\infty,X_{\ga,\mu}}$ can be estimated by \eqref{LWP6}. In summary, if we choose $r>0, T>0$ and $\ep>0$ sufficiently small, we obtain a constant $c=c(u_0)>0$ such that the estimate
    \begin{equation}\label{LWP7}
    ||\calT_{u_1}v_1-\calT_{u_2}v_2||_{\E_{1,\mu}(0,T)}\le\frac{1}{2}||v_1-v_2||_{\E_{1,\mu}(0,T)}+c|u_1-u_2|_{X_{\ga,\mu}},
    \end{equation}
is valid for all $u_1,u_2\in \bar{B}_{\ep}^{X_{\ga,\mu}}(u_0)$ and $v_1\in \B_{r,T,u_1}$, $v_2\in \B_{r,T,u_2}$. In the very special case $u_1=u_2$, \eqref{LWP7} yields the contraction mapping property of $\calT_{u_1}$ on $\B_{r,T,u_1}$. Now we are in a position to apply Banach's fixed point theorem to obtain a unique fixed point $\tilde{u}\in\B_{r,T,u_1}$ of $\calT_{u_1}$, i.e. $\calT_{u_1}\tilde{u}=\tilde{u}$. Therefore $\tilde{u}\in\B_{r,T,u_1}$ is the unique local solution to \eqref{LWP1}. Furthermore, if $u(t,u_1)$ and $u(t,u_2)$ denote the solutions of \eqref{LWP1} with initial values $u_1,u_2\in \bar{B}_{\ep}^{X_{\ga,\mu}}(u_0)$, respectively, the last assertion of the theorem follows from \eqref{LWP7}. The proof is complete.

\epr

\noindent
The next result provides information about the continuation of local solutions.
\begin{cor}\label{LWPcor1}
Let the assumptions of Theorem \ref{LWPthm} be satisfied and assume that $A(v)\in \calM\calR_{p}(X_1,X_0)$ for all $v\in V_\mu$. Then the solution $u(t)$ of \eqref{LWP1} has a maximal interval of existence $J(u_0)=[0,t^+(u_0))$.
\end{cor}
\bpr
Given $u_0\in X_{\gamma,\mu}$, Theorem \ref{LWPthm} yields some $T_1>0$ and a unique solution $\bar{u}\in \E_{1,\mu}(0,T_1)\cap C([0,T];V_\mu)$ of \eqref{LWP1}. Next, we apply Theorem \ref{LWPthm} to \eqref{LWP1} with initial value $\bar{u}(T_1)\in V_\mu$
to obtain some $T_2>0$ and a unique solution $\tilde{u}\in \E_{1,\mu}(0,T_2)\cap C([0,T_2];V_\mu)$. Let
    $$u(t):=\begin{cases}
            \bar{u}(t),&\quad t\in[0,T_1],\\
            \tilde{u}(t-T_1),&\quad t\in[T_1,T_1+T_2].
            \end{cases}
    $$
Then $u\in \E_{1,\mu}(0,T_1+T_2)\cap C([0,T_1+T_2];V_\mu)$, provided that
    \begin{equation}\label{LWP6a}
    \int_{T_1}^{T_1+T_2}|\tilde{u}(t-T_1)|_{1}^p\ t^{(1-\mu)p}\ dt+\int_{T_1}^{T_1+T_2}|\dot{\tilde{u}}(t-T_1)|_{0}^p\ t^{(1-\mu)p}\ dt<\infty,
    \end{equation}
since we already know $\bar{u}\in \E_{1,\mu}(0,T_1)$. To establish \eqref{LWP6a} it suffices to verify $\tilde{u}\in \E_{1}(0,T_2)$. Clearly, $\tilde{u}$ is a solution of the nonautonomous problem
    $$\dot{\tilde{u}}+\tilde{A}(t)\tilde{u}=\tilde{F}(t),\ t\in [0,T_2],\quad \tilde{u}(0)=\bar{u}(T_1),$$
where we have set $\tilde{A}(t):=A(\tilde{u}(t))$ and $\tilde{F}(t):=F(\tilde{u}(t))$. From \eqref{LWP2} it follows that
$\tilde{F}\in L_p(0,T_2;X_0)$ and $\tilde{A}\in C([0,T_2];\calB(X_1,X_0))$. The embedding
    $$\bar{u}\in \E_{1,\mu}(0,T_1)\hookrightarrow C((0,T_1];X_\gamma),$$
yields $\bar{u}(T_1)\in X_\gamma$. Therefore we may apply
\cite[Corollary 3.4]{JanBari} to obtain $\tilde{u}\in \E_{1}(0,T_2)$, whence $u\in \E_{1,\mu}(0,T_1+T_2)$ is the unique solution of \eqref{LWP1} on the interval $[0,T_1+T_2]$. Inductively this yields a maximal interval of existence $J(u_0):=[0,t^+(u_0))\subset J_0$, which is of course half sided open, since otherwise we could continue the solution beyond $t^+(u_0)$ with initial value $u(t^+(u_0))$.

\epr

\begin{rem}
Let $J=[0,T_0]$ a compact interval and denote by $\calM\calR_p(J;X_1,X_0)$ the class of all linear operators $A_0:X_1\to X_0$ such that for all $f\in L_p(J;X_0)$ there exists a unique solution $u\in H_p^1(J;X_0)\cap L_p(J;X_1)$ of
    $$\dot{u}+A_0u=f,\ t\in (0,T_0],\quad u(0)=0.$$
It is well-known that this properts does not depend on the length of the interval $J$, and that there exists a number $\kappa>0$ such that the implication
    $$A_0\in \calM\calR_p(J;X_1,X_0)\Rightarrow A_0+\kappa I\in \calM\calR_p(X_1,X_0)$$
holds, see e.g.\ Pr\"{u}ss \cite{JanBari}. In this sense the assumption $A(u_0)\in\calM\calR_p(X_1,X_0)$ in Theorem \ref{LWPthm} can be replaced by the somewhat weaker condition $A(u_0)\in\calM\calR_p(J;X_1,X_0)$, we simply have to add $\kappa u$ to both sides of \eqref{LWP1}.
\end{rem}

\section{Relative compactness of orbits}

Let $u_0\in V_\mu$ be given. Suppose that $(A,F)$ satisfy \eqref{LWP2} and $A(v)\in \calM\calR_p(J;X_1,X_0)$ for all $v\in V_\mu$ and for some $\mu\in (1/p,1)$, where $J=[0,T]$ or $J=\R_+$. In the sequel we assume that the unique solution of \eqref{LWP1} satisfies $u\in BC([\tau,t^+(u_0));V_\mu\cap X_{\gamma})$ for some $\tau\in (0,t^+(u_0))$ and
    \begin{equation}\label{GWPdist}
    \dist(u(t),\partial V_\mu)\ge\eta>0
    \end{equation}
for all $t\in J(u_0)$. Suppose furthermore that
    \begin{equation}\label{GWPemb}
    X_{\gamma}\compemb X_{\gamma,\mu},\quad \mu\in (1/p,1).
    \end{equation}
It follows from the boundedness of $u(t)$ in $X_\gamma$ that the set $\{u(t)\}_{t\in J(u_0)}\subset V_\mu$ is relatively compact in $X_{\gamma,\mu}$, provided $\mu\in (1/p,1)$. By \eqref{GWPdist} it holds that $\calV:=\overline{\{u(t)\}}_{t\in J(u_0)}$ is a real subset of $V_\mu$. Applying Theorem \ref{LWPthm} we find for each $v\in \calV$ numbers $\ep(v)>0$ and $\de(v)>0$ such that $B_{\ep(v)}^{X_{\gamma,\mu}}(v)\subset V_\mu$ and all solutions of \eqref{LWP1} which start in $B_{\ep(v)}^{X_{\gamma,\mu}}(v)$ have the common interval of existence $[0,\de(v)]$. Therefore the set
    $$\bigcup_{v\in \calV}B_{\ep(v)}^{X_{\gamma,\mu}}(v)$$
is an open covering of $\calV$ and by compactness of $\calV$ there exist $N\in\N$ and $v_k\in\calV$, $k=1,\ldots,N$, such that
    $$\calU:=\bigcup_{k=1}^N B_{\ep_k}^{X_{\gamma,\mu}}(v_k)\supset\calV=\overline{\{u(t)\}}_{t\in J(u_0)}\supset\{u(t)\}_{t\in J(u_0)},$$
where $\ep_k:=\ep(v_k)$, $k=1,\ldots,N$. To each of these balls corresponds an interval of existence $[0,\de_k]$, $\de_k>0$, $k=1,\ldots,N$.
Consider the problem
    \begin{equation}\label{GWP4}
    \dot{v}+A(v)v=F(v),\ s>0,\quad v(0)=u(t),
    \end{equation}
where $t\in J(u_0)$ is fixed and let $\de:=\min\{\de_k,\ k=1,\ldots,N\}$. Since $u(t)\subset \calU,\ t\in J(u_0)$, the solution of \eqref{GWP4} exists at least on the interval $[0,\de]$. By uniqueness it holds that $v(s)=u(t+s)$ if $t+s\in J(u_0)$, $t\in J(u_0)$, $s\in[0,\de]$, hence $\sup J(u_0)=+\infty$, i.e.\ the solution exists globally.

By continuous dependence on the initial data, the solution operator $G_1:\mathcal{U}\to \E_{1,\mu}(0,\de)$, which assigns to each initial value $u_1\in \mathcal{U}$ a unique solution $v(\cdot,u_1)\in \E_{1,\mu}(0,\de)$, is continuous. Furthermore
    $$(\de/2)^{1-\mu}||v||_{\E_{1}(\de/2,\de)}\le ||v||_{\E_{1,\mu}(\de/2,\de)}\le ||v||_{\E_{1,\mu}(0,\de)},\ \mu\in (1/p,1),$$
wherefore the mapping $G_2:\E_{1,\mu}(0,\de)\to \E_{1}(\de/2,\de)$ with $v\mapsto v$ is continuous. Finally
    $$|v(\delta)|_{X_{\gamma}}\le||v||_{BUC((\de/2,\de);X_{\gamma})}\le C(\delta)||v||_{\E_1(\de/2,\de)},$$
hence the mapping $G_3:\E_1(\de/2,\de)\to X_{\gamma}$ with $v\mapsto v(\de)$ is continuous. This yields the continuity of the composition $G=G_3\circ G_2\circ G_1:\mathcal{U}\to X_{\gamma}$, whence $G(\{u(t)\}_{t\ge 0})=\{u(t+\de)\}_{t\ge 0}$ is relatively compact in $X_{\gamma}$, since the continuous image of a relatively compact set is relatively compact. Since the solution has relatively compact range in $X_\gamma$, it is an easy consequence that the $\om$-limit set
    $$\om(u_0):=\left\{v\in V_\mu\cap X_\gamma:\ \exists\ t_n\nearrow\infty\ \mbox{s.t.}\ u(t_n;u_0)\to v\ \mbox{in}\ X_\gamma\right\}$$
is nonempty, connected and compact. We summarize the preceding considerations in the following
\begin{thm}\label{GWPthm2}
Let $p\in (1,\infty)$ and let $J=[0,T]$ or $J=\R_+$. Suppose that $A(v)\in\calM\calR_{p}(J;X_1,X_0)$ for all $v\in V_\mu$ and let \eqref{LWP2} as well as \eqref{GWPemb} hold for some $\mu\in (1/p,1)$. Assume furthermore that the solution $u(t)$ of \eqref{LWP1} satisfies 
    $$u\in BC([\tau,t^+(u_0));V_\mu\cap X_{\gamma})$$ 
for some $\tau\in (0,t^+(u_0))$ and
    $$\dist(u(t),\pa V_\mu)\ge \eta>0$$
for all $t\in J(u_0)$. Then the solution exists globally and for each $\de>0$, the orbit $\{u(t)\}_{t\ge \de}$ is relatively compact in $X_{\gamma}$. If in addition $u_0\in V_\mu\cap X_{\gamma}$, then $\{u(t)\}_{t\ge 0}$ is relatively compact in $X_{\gamma}$.
\end{thm}

\subsection{A second order quasilinear problem}

In this subsection we show how to apply Theorem \ref{GWPthm2} to a certain class of second order quasilinear parabolic equations. To be precise, we consider the problem
    \begin{align}
    \begin{split}\label{appl1}
    \pa_t u-a(u,\nabla u):\nabla^2 u&=f(u,\nabla u),\ t>0,\ x\in\Om,\\
    u&=0,\ t>0,\ x\in\pa\Om,\\
    u(0)&=u_0,\ x\in\Om
    \end{split}
    \end{align}
where $\Om\subset\R^n$ is a bounded domain with boundary $\pa\Om\in C^2$, $f\in C^{1}(\R\times\R^n;\R)$, $a\in C^{1}(\R\times\R^n;\R^{n\times n})$ and $a(u,v)$ is symmetric and positive definite for each $(u,v)\in \R\times\R^n$. If $A,B\in\R^{n\times n}$, then $A:B$ stands for
    $$A:B=\sum_{i,j=1}^n a_{ij}b_{ij}={\rm tr}(AB^{\textsf{T}}),$$
which defines the standard inner product in the space of matrices $\R^{n\times n}$. Let us first rewrite \eqref{appl1} in the form \eqref{LWP1}. To this end we set $X_0=L_p(\Om)$,
    $$X_1=\{u\in W_p^2(\Om):u|_{\pa\Om}=0\},$$
where $u|_{\pa\Om}$ has to be understood in the sense of traces. In this situation, we have for $\mu\in (1/p,1]$
    $$X_{\ga,\mu}=(X_0,X_1)_{\mu-1/p,p}=\begin{cases}
                                        \{u\in W_p^{2\mu-2/p}(\Om):u|_{\pa\Om}=0\},&\ \mbox{if}\ \mu p>3/2,\\
                                        W_p^{2\mu-2/p}(\Om),&\ \mbox{if}\ 1<\mu p<3/2,
                                        \end{cases}$$
see e.g.\ \cite{Gri69}. Let us assume that $p>n+2$, wherefore the embedding $W_p^{2-2/p}(\Om)\hookrightarrow C^1(\bar{\Om})$ is at our disposal. In this case there exists $\mu_0\in (1/p,1)$ such that
    $$W_p^{2-2/p}(\Om)\compemb W_p^{2\mu-2/p}(\Om)\hookrightarrow C^1(\bar{\Om}),\quad\mbox{if}\ \mu\in (\mu_0,1).$$
Indeed, the number $\mu_0\in (1/p,1)$ is given by
    $$\mu_0=\frac{1}{2}+\frac{n+2}{2p}=\frac{1}{p}+\frac{n+p}{2p},$$
provided $p>n+2$. For $\mu\in (\mu_0,1]$, we define $A:X_{\ga,\mu}\to \calB(X_0,X_1)$ and $F:X_{\ga,\mu}\to X_0$ by means of
    $$A(v)u(x):=a(v(x),\nabla v(x)):\nabla^2 u(x),\ x\in\Om,\ v\in X_{\ga,\mu},\ u\in X_1,$$
and
    $$F(v)(x):=f(v(x),\nabla v(x)),\ x\in\Om,\ v\in X_{\ga,\mu}.$$
From the regularity assumptions on $a$ and $f$ it follows that
    $$(A,F)\in C^{1-}(X_{\ga,\mu};\calB(X_1,X_0)\times X_0),\quad \mu\in (\mu_0,1].$$
Furthermore, by \cite{DHP1,DHP2}, we obtain $A(v)\in \calM\calR_p(J;X_1,X_0)$ for all $v\in X_{\gamma,\mu}$, $\mu\in (\mu_0,1]$, where $J=[0,T_0]$ is an arbitrary compact interval. By Theorem \ref{LWPthm} there exists a unique solution $u$ of \eqref{appl1} with maximal interval of existence $J(u_0)$, provided $u_0\in X_{\gamma,\mu}$. Assuming in addition $u\in BC(J(u_0);X_\gamma)$ we may apply Theorem \ref{GWPthm2} to the result
\begin{thm}
Let $n\in\N$, $p>n+2$, $\Omega\subset\R^n$ a bounded domain with boundary $\pa\Om\in C^2$ and let $u_0\in W_p^{2-2/p}(\Om)$ such that $u_0|_{\pa\Om}=0$. Assume in addition $f\in C^{1}(\R\times\R^n;\R)$ and $a\in C^{1}(\R\times\R^n;\R^{n\times n})$ with the property that $a(u,v)$ is symmetric and positive definite for each $(u,v)\in \R\times\R^n$. If the solution $u(t)$ of \eqref{appl1} satisfies
    $$u\in BC\left(J(u_0);W_p^{2-2/p}(\Om)\right),$$
then $u(t)$ exists globally, i.e.\ $J(u_0)=\R_+$ and the set $\{u(t)\}_{t\ge 0}$ is relatively compact in $W_p^{2-2/p}(\Om)$. Moreover, the $\om$-limit set
    $$\om(u_0):=\left\{v\in W_p^{2-2/p}(\Om):\ \exists\ t_n\nearrow\infty\ \mbox{s.t.}\ u(t_n;u_0)\to v\ \mbox{in}\ W_p^{2-2/p}(\Om)\right\}$$
is nonempty, connected and compact.
\end{thm}

\begin{rem}
For simplicity we supplied $\eqref{appl1}_1$ with a Dirichlet boundary condition. However, this boundary condition may be replaced by any other even nonlinear one (up to differential order one), as long as the Lopatinskii-Shapiro condition holds, which leads to maximal $L_p$-regularity (see e.g.\ \cite{DHP1,DHP2,LaPrSchn06}).
\end{rem}

\section{Long-Time Behavior}

In this section we investigate the long-time behavior of solutions to the quasilinear problem
    \begin{equation}\label{LTB1}
    \dot{u}+A(u)u=F(u),\ t>0,\quad u(0)=u_0,
    \end{equation}
where $(A,F)\in C^1(V_\mu;\calB(X_1,X_0)\times X_0)$ and $V_\mu\subset X_{\gamma,\mu}=(X_0,X_1)_{\mu-1/p,p},\ \mu\in (1/p,1)$ is open. We call $u_*$ an equilibrium of \eqref{LTB1} if $u_*\in V_\mu\cap X_1$ and $A(u_*)u_*=F(u_*)$. The following result has been proven in \cite[Theorem 6.1]{PSZ} in the classical setting, i.e. $\mu=1$.
\begin{thm}\label{thmPSZ}
Let $1<p<\infty$ and let $V_1\subset X_\gamma$ be open. Suppose $u_*\in V_1\cap X_1$ is an equilibrium of \eqref{LTB1}, and suppose that $(A,F)\in C^1(V_1;\calB(X_1,X_0)\times X_0)$. Suppose further that $A(u_*)$ has the property of maximal $L_p$-regularity. Let $A_0$ be the linearization of \eqref{LTB1} at $u_*$. Suppose that $u_*$ is \emph{normally hyperbolic}, i.e. assume that
    \begin{enumerate}
    \item near $u_*$ the set of equilibria $\calE\subset V_1\cap X_1$ is a $C^1$-manifold in $X_1$ of dimension $m\in\mathbb{N}_0$,
    \item the tangent space for $\calE$ at $u_*$ is given by $N(A_0)$,
    \item 0 is a semi-simple eigenvalue of $A_0$, i.e. $N(A_0)\oplus R(A_0)=X_0$,
    \item $\sigma(A_0)\cap i\R\subset\{0\}$, $\sigma(A_0)\cap \C_-\neq\emptyset$, $\sigma(A_0)\cap \C_+\neq\emptyset$.
    \end{enumerate}
Then for each sufficiently small $\rho>0$ there exists $\de\in(0,\rho]$ such that the unique solution $u(t)$ of \eqref{LTB1} with initial value $u_0\in B_\de^{X_\gamma}(u_*)$ either satisfies
    \begin{itemize}
    \item[(\textit{a})] $\dist_{X_\gamma}(u(t_0),\calE)>\rho$ for some finite time $t_0>0$, or
    \item[(\textit{b})] $u(t)$ exists on $\R_+$ and converges at an exponential rate to some $u_\infty\in\calE$ in $X_\gamma$ as $t\to\infty$.
    \end{itemize}
If $u_*$ is \emph{normally stable}, i.e.\ if in addition $\si(A_0)\cap\C_+=\emptyset$, then (a) does not occur.
\end{thm}
\begin{rem}
If $m=0$ then conditions (i)-(iv) of Theorem \ref{thmPSZ} imply $N(A_0)=\{0\}$, hence $R(A_0)=X_0$ and $\si(A_0)\cap i\R=\emptyset$, i.e.\ $u_*$ is \emph{hyperbolic}. The inverse function theorem then yields that $u_*$ is isolated in $V_1\cap X_1$, and so in this case Theorem \ref{thmPSZ} is contained in \cite[Theorem 7.1]{JanBari}.
\end{rem}
\noindent
It is our aim to extend this local result on qualitative behavior to a global one, under the slightly stronger assumption $(A,F)\in C^1(V_\mu;\calB(X_1,X_0)\times X_0)$ for some $\mu\in (1/p,1)$ and provided that \eqref{GWPemb} holds. Let $V_{\mu,\gamma}:=V_\mu\cap X_\gamma$.
Assume that $u\in BC(\R_+;V_{\mu,\gamma})$ is a global solution to \eqref{LTB1}, satisfying
    $$\dist(u(t),\partial V_{\mu})\ge\eta>0$$
for all $t\ge 0$. The mapping $(t,u_1)\mapsto S(t)u_1$, defined by $S(t)u_1=u(t,u_1),\ t\ge 0,\ u_1\in V_{\mu,\gamma}$ defines a semiflow in $V_{\mu,\gamma}$. Let $\Phi\in C(V_{\mu,\gamma};\R)$ be a strict Ljapunov function for $\{S(t)\}_{t\ge 0}$, that is
    \begin{itemize}\itemsep=0.15cm
    \item[($\Phi1$)] The function $t\mapsto \Phi(S(t)u_0)$ is nonincreasing, and
    \item[($\Phi2$)] If $\Phi(S(t)u_*)=\Phi(u_*)$ for all $t\ge 0$ then $u_*\in V_\mu\cap X_1$ is an equilibrium of \eqref{LTB1}.
    \end{itemize}
Theorem \ref{GWPthm2} yields that the orbit $\{u(t)\}_{t\ge 0}$ is relatively compact in $X_\gamma$. Hence the $\omega$-limit set
    \begin{equation}\label{LTB2}
    \omega(u_0)=\{v\in V_{\mu,\gamma}:\ \exists\ t_n\nearrow+\infty\ \mbox{s.t.}\ S(t_n)u_0\to v\ \mbox{in}\ X_\gamma,\ \mbox{as}\ n\to\infty\}
    \end{equation}
is nonempty, compact, and connected. Moreover, $\dist(S(t)u_0,\omega(u_0))\to 0$ in $X_\gamma$ as $t\to \infty$ and $\omega (u_0)\subset \calE\subset V_\mu\cap X_1$, by $(\Phi1)$ and $(\Phi2)$,  wherefore the set of equilibria is nonempty. Let $u_*\in\om(u_0)$. Then there exists a sequence $t_n\nearrow+\infty$ such that $S(t_n)u_0\to u_*$ in $X_\gamma$ as $t_n\to\infty$. Assuming that $u_*$ is normally hyperbolic and $t_n$ is large enough, Theorem \ref{thmPSZ} yields the convergence of $S(t)u_0$ to some equilibrium $u_\infty\in V_{\mu,\gamma}$ as $t\to\infty$. Uniqueness of the limit finally implies $u_\infty=u_*$. We obtain the following result.
\begin{thm}\label{LTBthm}
Let $p\in (1,\infty)$, $\mu\in (1/p,1)$, $V_\mu\subset X_{\gamma,\mu}$ be open, $V_{\mu,\gamma}=V_\mu\cap X_\gamma$ and assume that $(A,F)\in C^1(V_\mu;\calB(X_1,X_0)\times X_0)$ and \eqref{GWPemb} hold for some $\mu\in (1/p,1)$. Suppose furthermore that $u\in BC(\R_+;V_{\mu,\gamma})$ is a global solution to \eqref{LTB1}, satisfying
    $$\dist(u(t),\partial V_{\mu})\ge\eta>0$$
for all $t\ge 0$ and let $\Phi\in C(V_{\mu,\gamma};\R)$ be a strict Ljapunov function for \eqref{LTB1}. Then the $\om$-limit set, defined by \eqref{LTB2}, is nonempty, compact and connected. If in addition there exists $u_*\in\om(u_0)$ which is normally hyperbolic, then $\lim_{t\to\infty}u(t)=u_*$ in $X_\gamma$, $u_*\in V_{\mu}\cap X_1$ and $A(u_*)u_*=F(u_*)$.
\end{thm}

\section{The Mullins-Sekerka problem}

Let $\Omega\subset\R^n,\ n\ge 2$ be a bounded domain with a smooth boundary $\pa\Om$. Let $\Gamma_0\subset\Om$ be a compact connected hypersurface in $\Om$ which divides $\Om$ into two disjoint sets $\Om_1^0$ (liquid phase) and $\Om_2^0$ (solid phase) such that $\pa\Om_1^0=\Gamma_0$ and $\pa\Om_2^0=\Gamma_0\cup \pa\Om$. We regard $\Gamma_0$ as the initial state of a time dependent family of hypersurfaces $\{\Gamma(t)\}_{t\ge 0}$ and denote by $\Gamma(t)$ its position at time $t>0$. Let $V(t,\cdot)$ and $\kappa(t,\cdot)$ be the normal velocity and the mean curvature of $\Gamma(t)$, and let $\Om_1(t)$ and $\Om_2(t)$ be the two disjoint regions in $\Om$ which are separated by $\Gamma(t)$, such that $\pa\Om_1(t)=\Gamma(t)$ and $\partial\Om_2(t)=\Gamma(t)\cup\pa\Om$. Let further $\nu_\Gamma(t,\cdot)$ be the outer unit normal field on $\Gamma(t)$ w.r.t.\ $\Om_1(t)$ and let $\nu(\cdot)$ be the outer unit normal field on $\pa\Om$. The two-phase Mullins-Sekerka problem consists in finding a family $\{\Gamma(t)\}_{t\ge 0}$ of hypersurfaces satisfying
    \begin{equation}\label{MS1}
    V=[\![\partial_{\nu_\Gamma} u_\kappa]\!],\quad t>0,\quad \Gamma(0)=\Gamma_0,
    \end{equation}
where $u_\kappa=u_\kappa(t,\cdot)$ is, for each fixed $t\ge 0$, the unique solution of the elliptic boundary value problem
    \begin{align}
    \begin{split}\label{MS2}
    \De u&=0,\quad x\in\Om\backslash\Gamma(t),\\
    u&=\kappa,\quad x\in\Gamma(t),\\
    \pa_\nu u&=0,\quad x\in\pa\Om.
    \end{split}
    \end{align}
Here $[\![\partial_{\nu_\Gamma} u_\kappa]\!]:=\pa_{\nu_\Gamma}u_\kappa^2-\pa_{\nu_\Gamma}u_\kappa^1$ stands for the jump of the normal derivative of $u_\kappa$ across the interface $\Gamma(t)$. In order to reformulate the Mullins-Sekerka problem as a quasilinear evolution equation in an abstract $L_p$-setting, we need some preliminaries from differential geometry. Let $\Sigma\subset\Omega$ be a real analytic ($C^\om$-) hypersurface such that $\Sigma$ divides $\Omega$ into two disjoint, open, connected sets $\Omega_1$ and $\Omega_2$, the interior and the exterior of $\Sigma$. It is well-known that $\Sigma$ admits a \emph{tubular} neighborhood, which means that there is a number $a>0$ such that the map
\begin{equation*}
\Lambda: \Sigma \times (-a,a)\to \R^n,\quad \Lambda(p,r):= p+r\nu_\Sigma(p),
\end{equation*}
is a $C^\om$-diffeomorphism from $\Sigma \times (-a,a)$
onto its image $U_a:=R(\Lambda)$. The inverse
$$\Lambda^{-1}:R(\Lambda)\mapsto \Sigma\times (-a,a)$$ of this map
is conveniently decomposed as
$$\Lambda^{-1}(x)=(\Pi(x),d_\Sigma(x)),\quad x\in R(\Lambda).$$
Here $\Pi(x)$ means the orthogonal projection of $x$ to $\Sigma$ and $d_\Sigma(x)$ the signed
distance from $x$ to $\Sigma$; so $|d_\Sigma(x)|={\rm dist}(x,\Sigma)$ and $d_\Sigma(x)<0$ if and only if
$x\in \Omega_1$. In particular we have $R(\Lambda)=\{x\in \R^n:\, {\rm dist}(x,\Sigma)<a\}$.
Note that one the one side an upper bound for $a$ is determined by the curvatures of $\Sigma$, i.e.\ we must have
$$0<a<\min\{1/\kappa_j(p): j=1,\ldots,n-1,\; p\in\Sigma\},$$
where $\kappa_j(p)$ mean the principal curvatures of $\Sigma$ at $p\in\Sigma$. On the other side, $a$ is also connected to the topology of $\Sigma$,
which can be expressed as follows. Since $\Sigma$ is a compact manifold of
dimension $n-1$ it satisfies the ball condition, which means that
there is a radius $r_\Sigma>0$ such that for each point $p\in \Sigma$
there are $x_j\in \Omega_j$, $j=1,2$, such that $B_{r_\Sigma}(x_j)\subset \Omega_j$, $j=1,2$, and
$\bar{B}_{r_\Sigma}(x_j)\cap\Sigma=\{p\}$. Choosing $r_\Sigma$ maximal,
we then must also have $a<r_\Sigma$.

In case $\Gamma_0\subset R(\Lambda)$, we may use the map $\Lambda$ to parameterize the unknown free
boundary $\Gamma(t)$ over $\Sigma$ by means of a height function $\rho(t,p)$ with $|\rho|_\infty<a$ via
$$\Gamma(t): p\mapsto p+ \rho(t,p)\nu_\Sigma(p),\quad p\in\Sigma,\; t\geq0,$$
for small $t\geq0$, at least. We extend this diffeomorphism to all of $\bar{\Omega}$ by means of
$$ \Theta_h(t,x) = x +\chi(d_\Sigma(x)/a)\rho(t,\Pi(x))\nu_\Sigma(\Pi(x)).$$
Here $\chi$ denotes a suitable cut-off function; more precisely, $\chi\in C^\infty(\R)$,
$0\leq\chi\leq 1$, $\chi(s)=1$ for $|s|<1/3$, and $\chi(r)=0$ for $|s|>2/3$.
This way $\Omega\setminus\Gamma(t)$ is transformed to the fixed domain $\Omega\setminus\Sigma$. This is known as the Hanzawa transform. Following \cite{ES98} we obtain for the transformed problem \eqref{MS1} the initial value problem
    \begin{equation}\label{MS3}
    \dot{\rho}+B(\rho)S(\rho)K(\rho)=0,\quad t>0,\quad \rho(0)=\rho_0,
    \end{equation}
on $\Sigma$. Here $S(\rho)g$ is the solution of the transformed elliptic boundary value problem
    \begin{align}
    \begin{split}\label{MS4}
    \calA(\rho) v&=0,\quad x\in\Om\backslash\Sigma,\\
    v&=g,\quad x\in\Sigma,\\
    \pa_\nu v&=0,\quad x\in\pa\Om,
    \end{split}
    \end{align}
where $\calA(\rho)$ means the transformed Laplacian and $K(\rho)$ resp.\ $B(\rho)$ denote the transformed mean curvature operator resp.\ the transformed jump of the normal derivative. We want to study \eqref{MS3} in an $L_p$-setting. Let $p>(n+3)/2$ and define
    $$X_0=W_p^{1-1/p}(\Sigma),\quad X_1=W_p^{4-1/p}(\Sigma).$$
We consider \eqref{MS3} as an evolution equation in the space $L_{p,\mu}(J;X_0)$, where $J=[0,T]$, $T>0$ and $\mu\in (1/p,1]$. To be precise, we are looking for solutions in the maximal regularity class
    $$H_{p,\mu}^1(J;X_0)\cap L_{p,\mu}(J;X_1).$$
The corresponding (weighted) trace space is given by the real interpolation method and reads
    $$X_{\gamma,\mu}=(X_0,X_1)_{\mu-1/p,p}=W_p^{3\mu+1-4/p}(\Sigma).$$
Since $p>(n+3)/2$, the Sobolev embedding
    $$X_{\gamma}\compemb X_{\gamma,\mu}\hookrightarrow C^2(\Sigma)$$
holds for $\mu\in (\mu_0,1)$ with a sufficiently large $\mu_0\in (1/p,1)$. Here the number $\mu_0\in (1/p,1)$ is given by
    $$\mu_0=\frac{1}{3}+\frac{n+3}{3p}=\frac{1}{p}+\frac{n+p}{3p},$$
provided $p>(n+3)/2$. Note that we can choose the real analytic hypersurface $\Sigma$ in such a way that $|\rho_0|_{X_{\gamma,\mu}}\le\de$ with a sufficiently small $\de>0$. Therefore we define the set $V_\mu$ from Theorem \ref{LWPthm} to be the open ball $B_\de^{X_{\ga,\mu}}(0)\subset X_{\gamma,\mu},\ \mu\in (\mu_0,1]$. It is well known that $K(\rho)$ can be decomposed as
    $$K(\rho)=P(\rho)\rho+Q(\rho),$$
where $P\in C^1(V_\mu;\calB(X_1;Y)), Y:=W_p^{2-1/p}(\Sigma)$, is a differential operator of second order and $Q\in C^1(V_\mu;Y)$ contains only first order terms. Moreover, the (transformed) two-phase Dirichlet-to-Neumann operator $B(\rho)S(\rho)$ has the property $BS\in C^1(V_\mu;\calB(Y,X_0))$. This yields
    $$A:=BSP\in C^1(V_\mu;\calB(X_1,X_0)),$$
and $F:=BSQ\in C^1(V_\mu;X_0)$. Now we take care about the maximal regularity property of $A(\rho_0)$, $\rho_0\in V_\mu$. In other words we want to show that for $J=[0,T]$, $T>0$, and any $f\in L_p(J;X_0)$ the problem
    $$\dot{\sigma}+A(\rho_0)\sigma=f,\quad t>0,\quad \sigma(0)=0,$$
admits a unique solution $\si\in H_p^1(J;X_0)\cap L_p(J;X_1)$. For this purpose we show first that $A(0)$ has this property. Since $P(0)=-\De_{\Sigma}$, we have $A(0)\rho=[\![\partial_{\nu_{\Sigma}}S(\De_{\Sigma}\rho)]\!]$, where $Sg$ is the solution operator of
    \begin{align}\label{MS5}
    \begin{split}
    \De v&=0,\quad x\in\Om\backslash\Sigma,\\
    v&=g,\quad x\in\Sigma,\\
    \pa_\nu v&=0,\quad x\in\pa\Om.
    \end{split}
    \end{align}
Hence the maximal regularity property of $A(0)$ follows from \cite[Section 4]{PSZ}. Since $A(\rho_0)=A(0)+A(\rho_0)-A(0)$, maximal $L_p$-regularity of $A(\rho_0)$ follows by a perturbation argument, provided $|\rho_0|_{X_{\gamma,\mu}}$ is small enough. By Theorem \ref{LWPthm} there exist $T=T(\rho_0)>0$ and $\ep=\ep(\rho_0)>0$, such that $\bar{B}_{\ep}^{X_{\gamma,\mu}}(\rho_0)\subset V_\mu$ and such that the problem
    $$\dot{\rho}+B(\rho)S(\rho)K(\rho)=0,\quad t>0,\quad \rho(0)=\rho_1$$
has a unique solution
    $$\rho(\cdot,\rho_1)\in H_{p,\mu}^1(0,T;X_0)\cap L_{p,\mu}(0,T;X_1)\cap C([0,T];V_\mu),$$
on $[0,T]$, for any initial value $\rho_1\in \bar{B}_\ep^{X_{\gamma,\mu}}(\rho_0)$. Furthermore there exists a constant $c=c(\rho_0)>0$ such that for all $\rho_1,\rho_2\in \bar{B}_\ep^{X_{\gamma,\mu}}(\rho_0)$ the estimate
    $$||\rho(\cdot,\rho_1)-\rho(\cdot,\rho_1)||_{\E_{1,\mu}(0,T)}\le c|\rho_1-\rho_2|_{X_{\gamma,\mu}}$$
is valid. By regularization we even have $\rho(t;\rho_1)\in V_\mu\cap X_{\gamma}$ for all $t\in (0,T]$. The result on local well-posedness of \eqref{MS1} reads as follows.
\begin{thm}\label{MSthm1}
Let $(n+3)/2<p<\infty$ and $\mu\in(\mu_0,1]$. For each $\Gamma_0\in W_p^{3\mu+1-4/p}$, the Mullins-Sekerka problem \eqref{MS1} has a unique solution $\Gamma(t)\in W_p^{4-4/p}$ on a possibly small time interval $(0,T]$. The solution depends continuously on the initial data.
\end{thm}
By a proper choice of the real analytic hypersurface $\Sigma$, reparametrization and successive application of Theorem \ref{MSthm1} yields a maximal interval of existence $J(\Gamma_0)=[0,t^+(\Gamma_0))$ for the solution $\Gamma(t)$ of $\eqref{MS1}$. In order to investigate global existence in time as well as long-time behaviour, we need some more facts from differential geometry. First of all, we recall that the set of all $C^2$-hypersurfaces which are contained in $\Omega$, form a $C^2$-manifold, which we denote by $\calM H^2(\Omega)$. A metric on $\calM H^2(\Omega)$ can be defined as follows. The Hausdorff metric $d_H$, defined on the set $\calK$ of all compact subsets of $\R^n$ is given by
    $$d_H(K_1,K_2)=\max\left\{\sup_{x\in K_1}d(x,K_2),\sup_{x\in K_2}d(K_1,x)\right\},\quad K_1,K_2\in\calK.$$
For $\Sigma_1,\Sigma_2\in\calM H^2(\Omega)$, we define
    $$d(\Sigma_1,\Sigma_2):=d_H(\calN^2\Sigma_1,\calN^2\Sigma_2),$$
as a metric on $\calM H^2(\Omega)$, where $\calN^2\Sigma$ stands for the second normal bundle of the hypersurface $\Sigma\in\calM H^2(\Omega)$, which is given by
    $$\calN^2\Sigma=\left\{p,\nu_\Sigma(p),\nabla_\Sigma\nu_\Sigma(p):p\in\Sigma\right\}.$$
The charts for $\calM H^2(\Omega)$ are the parameterizations over real analytic hypersurfaces $\Sigma\subset\Omega$. In this sense $\calM H^2(\Omega)$ becomes a Banach manifold.

Next, we want to show that each $\Sigma\in \calM H^2(\Omega)$ has a \emph{level function}. Let $U_a$ be the tube for $\Sigma$ and assume w.l.o.g.\ that $a\le 1$. We may then define a function $\varphi_\Sigma\in C^2(\R^n)$ by means of
$$\varphi_\Sigma(x) = g(d_\Sigma(x)),\quad x\in\R^n,$$
where
$$g(s)=s\chi(s/a)+ (1-\chi(s/a)){\rm sgn}\, s,\quad s\in \R,$$
and $\chi\in C^\infty$ is defined as above.
Then it is easy to see that $\Sigma=\varphi_\Sigma^{-1}(0)$, and $\nabla \varphi_\Sigma(x)=\nu_\Sigma(x)$, for each $x\in \Sigma$.

Consider the subset $\calM H^2(\Omega,r)$ of $\calM H^2(\Omega)$ which consists of all $\Gamma\in\calM H^2(\Omega)$ such that $\Gamma\subset\Omega$ satisfies the ball condition with radius $r>0$. This implies in particular ${\rm dist}(\Gamma,\partial\Omega)\geq r$  and all principal curvatures of $\Gamma\in\calM H^2(\Omega,r)$ are bounded by $r$. Further, the level functions $\varphi_\Gamma = g\circ d_{\Gamma}$ are well defined for $\Gamma\in\calM H^2(\Omega,r)$, and form a bounded subset of $C^2(\bar{\Omega})$. The map $\Phi:\calM H^2(\Omega,r)\to C^2(\bar{\Omega})$ defined by $\Phi(\Gamma)=\varphi_\Gamma$ is an isomorphism of the metric space $\calM H^2(\Omega,r)$ onto $\Phi(\calM H^2(\Omega,r))\subset C^2(\bar{\Omega})$.
Let $s-(n-1)/p>2$; for $\Gamma_j\in\calM H^2(\Omega,r)$, $j=1,2$, we define $\Gamma_j\in W^s_p(\Omega,r)$ if $\varphi_{\Gamma_j}\in W^s_p(\Omega)$ and $\dist_{W_p^s}(\Gamma_1,\Gamma_2):=|\ph_{\Gamma_1}-\ph_{\Gamma_2}|_{W_p^s(\Omega)}$. In this case the local charts for $\Gamma\in\calM H^2(\Omega,r)$ can be chosen of class $W^s_p$ as well. Finally, a subset $K\subset W^s_p(\Omega,r)$ is said to be (relatively) compact, if $\Phi(K)\subset W^s_p(\Omega)$ is (relatively) compact.

With the help of the preceding considerations we may define an appropriate phase-manifold $\calP M$ for the two-phase Mullins-Sekerka problem
by means of
\begin{equation*}
\calP M:=\{\Gamma\in \calM H^2(\Omega):\Gamma\in W^{4-4/p}_p\}.
\end{equation*}
It is an easy consequence of Theorem \ref{MSthm1} that the solution $\Gamma(t)$ of \eqref{MS1} defines a local semiflow in $\calP M$ on the maximal interval $J(\Gamma_0)$.

Let us next discuss the equilibria of \eqref{MS1}. To this end we define a functional $\phi$ by means of
    \begin{equation}\label{MSLjap}
    \phi(\Gamma(t))=\int_{\Gamma(t)}d\Gamma={\rm mes}\Gamma(t).
    \end{equation}
Then the time derivative of $\phi(\Gamma(t))$ reads
    \begin{align*}
    \frac{d}{dt}\phi(\Gamma(t))&=-\int_{\Gamma(t)}V(t,x)\kappa(t,x)d\Gamma=-\int_{\Gamma(t)}[\![\partial_{\nu_\Gamma} u(t,x)]\!]u(t,x)d\Gamma\\
    &=-\int_{\Omega(t)}\diver(\nabla u u)dx=-\int_{\Omega(t)}|\nabla u|^2dx\le 0,
    \end{align*}
where me made use of the transport theorem and \eqref{MS1},\eqref{MS2}. This shows that $\phi$ is a Ljapunov functional for \eqref{MS1} and it is even a strict one, since $\dot{\phi}(\Gamma(t))=0$ if and only if $u=\kappa$ is constant, hence $V=0$. Since $\Om$ is bounded, it follows that $\Gamma$ must be a sphere $S_R(x_0)\subset\Omega$ with radius $R>0$ and center $x_0\in\Omega$. If conversely $V=0$, then $\kappa$ is constant. In other words, the set of equilibria $\calE$ of the Mullins-Sekerka problem \eqref{MS1} is given by
    $$\calE=\{S_R(x_0):R>0,\ \bar{B}_R(x_0)\subset\Omega\}.$$
Basically there are two facts which prevent the solution from existence on $\R_+$, namely
    \begin{itemize}
    \item \textbf{Regularity}: the norm of $\Gamma(t)$ in $W_p^{4-4/p}$ becomes unbounded as $t\nearrow t^+(\Gamma_0)$;
    \item \textbf{Geometry}: the topology of the interface $\Gamma(t)$ changes,
    or the interface touches the boundary of $\Omega$.
    \end{itemize}
We say that the solution $\Gamma(t)$ satisfies a \emph{uniform ball condition},
if there is a radius $r>0$ such that $\Gamma(J(\Gamma_0))\subset \calM H^2(\Omega,r)$.
Note that this condition bounds the curvature of $\Gamma(t)$, and prevents it to touch the outer
boundary $\partial \Omega$, or to undergo topological changes. The main result of this section reads as follows.
\begin{thm}\label{MSthm2}
Let $(n+3)/2<p<\infty$ and let $\Gamma(t)$ be a solution of the Mullins-Sekerka problem \eqref{MS1} on the maximal time interval $J(\Gamma_0)=[0,t^+(\Gamma_0))$. Assume furthermore that
    \begin{enumerate}
    \item $|\Gamma(t)|_{W_p^{4-4/p}}\le M<\infty$ for all $t\in J(\Gamma_0)$, and
    \item $\Gamma(t)$ satisfies a uniform ball condition for all $t\in J(\Gamma_0)$.
    \end{enumerate}
Then $J(\Gamma_0)=\R_+$, i.e.\ the solution exists globally, and $\Gamma(t)$ converges in $\calP M$ to an equilibrium $\Gamma_\infty\in\calE$ at an exponential rate. To be precise, there exists $\omega>0$ such that
    $$e^{\omega t}{\rm dist}_{W_p^{4-4/p}}(\Gamma(t),\Gamma_\infty)\to 0$$
as $t\to\infty$.
\end{thm}
\bpr
Assume that (i) and (ii) are valid. Then $\Gamma(J(\Gamma_0))\subset W^{4-4/p}_p(\Omega,r)$ is bounded, hence relatively compact in
$W^{3\mu+1-4/p}_p(\Omega,r)$ for $\mu\in (\mu_0,1)$. Thus we may cover this set by finitely many balls with centers $\Sigma_k$ which are real analytic such that
${\rm dist}_{W_p^{3\mu+1-4/p}}(\Gamma(t),\Sigma_j)\leq \delta$ for some $j=j(t)$, $t\in J(\Gamma_0)$, $\mu\in (\mu_0,1)$. Let $J_k=\{t\in J(\Gamma_0):\, j(t)=k\}$. Using for each $k$ a Hanzawa-transformation, we may employ Theorem \ref{MSthm1} to obtain solutions
$\Gamma^1$ with initial configurations $\Gamma(t)$ in the phase manifold $\calP M$ on a common time interval say $[0,T]$, and by uniqueness we have
$\Gamma^1(t)=\Gamma(t+T)$, $t+T\in J(\Gamma_0)$. Since the solution depends continuously on the initial data, the set $\Gamma(J(\Gamma_0))$ is relatively compact in $\calP M$. In particular this yields $J(\Gamma_0)=\R_+$ and the orbit $\Gamma(\R_+)$ is relatively compact in $\calP M$. As we already know, the mapping $\phi$ defined by \eqref{MSLjap} is a strict Ljapunov functional, hence the limit set $\omega(\Gamma_0)$
of a solution is contained in the set $\calE$ of equilibria. By compactness $\omega(\Gamma_0)\subset \calP M$ is non-empty, hence the solution comes close to $\calE$, i.e.\ there is a sequence $t_n\to\infty$ such that $\Gamma(t_n)\to \Gamma_\infty\in\calE$. For sufficiently large $t_n$ we parameterize $\Gamma(t_n)$ over $\Gamma_\infty$ by a height function $\rho(t_n;\cdot)$. By \cite[Section 4]{PSZ} all conditions of Theorem \ref{LTBthm} are satisfied for the corresponding transformed equation \eqref{MS3}. Therefore an application of Theorem \ref{LTBthm} completes the proof.
\epr

\begin{rem}
The conditions (i) and (ii) of Theorem \ref{MSthm2} are also necessary for global existence and convergence of $\Gamma(t)$ to some sphere $\Gamma_\infty\in\calE$. This follows from a compactness argument.
\end{rem}

\bibliographystyle{amsplain}
\bibliography{WilkeBib2}

\providecommand{\bysame}{\leavevmode\hbox to3em{\hrulefill}\thinspace}
\providecommand{\MR}{\relax\ifhmode\unskip\space\fi MR }
\providecommand{\MRhref}[2]{%
  \href{http://www.ams.org/mathscinet-getitem?mr=#1}{#2}
}
\providecommand{\href}[2]{#2}
\begin{thebibliography}{10}

\bibitem{Ama88}
H.~Amann, \emph{Dynamic theory of quasilinear parabolic equations. {I}.
  {A}bstract evolution equations}, Nonlinear Anal. \textbf{12} (1988), no.~9,
  895--919.

\bibitem{Ama89}
\bysame, \emph{Dynamic theory of quasilinear parabolic systems. {III}. {G}lobal
  existence}, Math. Z. \textbf{202} (1989), no.~2, 219--250.

\bibitem{Ama90}
\bysame, \emph{Dynamic theory of quasilinear parabolic equations. {II}.
  {R}eaction-diffusion systems}, Differential Integral Equations \textbf{3}
  (1990), no.~1, 13--75.

\bibitem{Ama93}
\bysame, \emph{Nonhomogeneous linear and quasilinear elliptic and parabolic
  boundary value problems}, Function spaces, differential operators and
  nonlinear analysis ({F}riedrichroda, 1992), Teubner-Texte Math., vol. 133,
  Teubner, Stuttgart, 1993, pp.~9--126.

\bibitem{Ama05}
\bysame, \emph{Quasilinear parabolic problems via maximal regularity}, Adv.
  Differential Equations \textbf{10} (2005), no.~10, 1081--1110.

\bibitem{Ang90}
S.~B. Angenent, \emph{Nonlinear analytic semiflows}, Proc. Roy. Soc. Edinburgh
  Sect. A \textbf{115} (1990), no.~1-2, 91--107.

\bibitem{BrHuLu00}
C.-M. Brauner, J.~Hulshof, and A.~Lunardi, \emph{A general approach to
  stability in free boundary problems}, J. Differential Equations \textbf{164}
  (2000), no.~1, 16--48.

\bibitem{ChFaSch09}
R.~Chill, E.~Fa{\v{s}}angov{\'a}, and R.~Sch{\"a}tzle, \emph{Willmore blowups
  are never compact}, Duke Math. J. \textbf{147} (2009), no.~2, 345--376.

\bibitem{CleLi93}
Ph. Cl{\'e}ment and S.~Li, \emph{Abstract parabolic quasilinear equations and
  application to a groundwater flow problem}, Adv. Math. Sci. Appl. \textbf{3}
  (1993/94), no.~Special Issue, 17--32.

\bibitem{CleSi01}
Ph. Cl{\'e}ment and G.~Simonett, \emph{Maximal regularity in continuous
  interpolation spaces and quasilinear parabolic equations}, J. Evol. Equ.
  \textbf{1} (2001), no.~1, 39--67.

\bibitem{DHP1}
R.~Denk, M.~Hieber, and J.~Pr{\"u}ss, \emph{{$\calR$}-boundedness, {F}ourier
  multipliers and problems of elliptic and parabolic type}, Mem. Amer. Math.
  Soc. \textbf{166} (2003), no.~788, viii+114.

\bibitem{DHP2}
\bysame, \emph{Optimal {$L\sb p$}-{$L\sb q$}-estimates for parabolic boundary
  value problems with inhomogeneous data}, Math. Z. \textbf{257} (2007), no.~1,
  193--224.

\bibitem{Esch94}
J.~Escher, \emph{On quasilinear fully parabolic boundary value problems},
  Differential Integral Equations \textbf{7} (1994), no.~5-6, 1325--1343.

\bibitem{ES98}
J.~Escher and G.~Simonett, \emph{A center manifold analysis for the
  {M}ullins-{S}ekerka model}, J. Differential Equations \textbf{143} (1998),
  no.~2, 267--292.

\bibitem{Gri69}
P.~Grisvard, \emph{\'{E}quations diff\'erentielles abstraites}, Ann. Sci.
  \'Ecole Norm. Sup. (4) \textbf{2} (1969), 311--395.

\bibitem{LSU}
O.~A. Lady{\v{z}}enskaja, V.~A. Solonnikov, and N.~N. Ural'ceva, \emph{Linear
  and quasilinear equations of parabolic type}, Translated from the Russian by
  S. Smith. Translations of Mathematical Monographs, Vol. 23, American
  Mathematical Society, Providence, R.I., 1967.

\bibitem{LaPrSchn06}
Y.~Latushkin, J.~Pr{\"u}ss, and R.~Schnaubelt, \emph{Stable and unstable
  manifolds for quasilinear parabolic systems with fully nonlinear boundary
  conditions}, J. Evol. Equ. \textbf{6} (2006), no.~4, 537--576.

\bibitem{LaPrSchn08}
\bysame, \emph{Center manifolds and dynamics near equilibria of quasilinear
  parabolic systems with fully nonlinear boundary conditions}, Discrete Contin.
  Dyn. Syst. Ser. B \textbf{9} (2008), no.~3-4, 595--633.

\bibitem{Lun95}
A.~Lunardi, \emph{Analytic semigroups and optimal regularity in parabolic
  problems}, Progress in Nonlinear Differential Equations and their
  Applications, 16, Birkh\"auser Verlag, Basel, 1995.

\bibitem{JanBari}
J.~Pr{\"u}ss, \emph{Maximal regularity for evolution equations in {$L\sb
  p$}-spaces}, Conf. Semin. Mat. Univ. Bari \textbf{285} (2002), 1--39 (2003).

\bibitem{PrSi04}
J.~Pr{\"u}ss and G.~Simonett, \emph{Maximal regularity for evolution equations
  in weighted {$L\sb p$}-spaces}, Arch. Math. (Basel) \textbf{82} (2004),
  no.~5, 415--431.

\bibitem{PSZ}
J.~Pr{\"u}ss, G.~Simonett, and R.~Zacher, \emph{On convergence of solutions to
  equilibria for quasilinear parabolic problems}, J. Differential Equations
  \textbf{246} (2009), 3902--3931.

\bibitem{Sim94}
G.~Simonett, \emph{Quasilinear parabolic equations and semiflows}, Evolution
  equations, control theory, and biomathematics ({H}an sur {L}esse, 1991),
  Lecture Notes in Pure and Appl. Math., vol. 155, Dekker, New York, 1994,
  pp.~523--536.

\bibitem{Sim95}
\bysame, \emph{Center manifolds for quasilinear reaction-diffusion systems},
  Differential Integral Equations \textbf{8} (1995), no.~4, 753--796.

\bibitem{Yag91}
A.~Yagi, \emph{Abstract quasilinear evolution equations of parabolic type in
  {B}anach spaces}, Boll. Un. Mat. Ital. B (7) \textbf{5} (1991), no.~2,
  341--368.

\end{thebibliography}

-----------------------------99614912995--
\end{document}